\documentclass[11pt]{article} 
\usepackage{amsfonts,amsmath,latexsym,amssymb,mathrsfs,amsthm,comment}
\usepackage{caption}
\usepackage{graphicx}
\usepackage{xcolor}

\let\OLDthebibliography\thebibliography
\renewcommand\thebibliography[1]{
  \OLDthebibliography{#1}
  \setlength{\parskip}{1pt}
  \setlength{\itemsep}{0pt plus 0.0ex}
}

\def\numberlikeadb{\global\def\theequation{\thesection.\arabic{equation}}}
\numberlikeadb
\newtheorem{theorem}{Theorem}[section]

\newtheorem{proposition}[theorem]{Proposition}
\newtheorem{remark}[theorem]{Remark}

\usepackage{color}

\usepackage{lscape}
\usepackage{caption}
\usepackage{multirow}
\allowdisplaybreaks
\begin{document}

\title{A Stein characterisation of the distribution of the product of correlated normal random variables
}
\author{Robert E. Gaunt\footnote{Department of Mathematics, The University of Manchester, Oxford Road, Manchester M13 9PL, UK, robert.gaunt@manchester.ac.uk; siqi.li-8@postgrad.manchester.ac.uk; heather.sutcliffe@manchester.ac.uk},\: Siqi L$\mathrm{i}^{*}$  and Heather L. Sutcliff$\mathrm{e}^{*}$}

\date{} 
\maketitle

\vspace{-7mm}

\begin{abstract} We obtain a Stein characterisation of the distribution of the product of two correlated normal random variables with non-zero means, and more generally the distribution of the sum of independent copies of such random variables. Our Stein characterisation is shown to naturally generalise a number of other Stein characterisations in the literature. From our Stein characterisation we derive recursive formulas for the moments of the product of two correlated normal random variables, and more generally the sum of independent copies of such random variables, which allows for efficient computation of higher order moments. 
\end{abstract}

\noindent{{\bf{Keywords:}}} Product of correlated normal random variables; Stein characterisation; Stein's method; moments

\noindent{{{\bf{AMS 2010 Subject Classification:}}} Primary 60E05; 62E15

\section{Introduction}

Let $(X, Y)$ be a bivariate normal random vector with mean vector $(\mu_X,\mu_Y)$, variances $(\sigma_X^2,\sigma_Y^2)$ and correlation coefficient $\rho$. Since the work of \cite{craig,wb32} in the 1930's, the distribution of the product $Z=XY$, which we denote by $PN(\mu_X,\mu_Y;\sigma_X^2,\sigma_Y^2;\rho)$, has received much attention in the statistics literature (see \cite{gaunt22,np16} for an overview of some of the literature), and has found numerous applications, with recent examples including condensed matter physics \cite{ach}, astrophysics \cite{cac} and chemical physics \cite{hey}. The mean $\overline{Z}_n=n^{-1}\sum_{i=1}^nZ_i$, where $Z_1,\ldots,Z_n$ are independent copies of $Z$, has also found applications in fields such as quantum cosmology \cite{gr}, electrical engineering \cite{ware} and astrophysics \cite{man}.

Recently, \cite{cui} derived the following formula for the probability density function (PDF) of $Z$: 
\begin{align}p_{Z}(x)&=\exp\bigg\{-\frac{1}{2(1-\rho^2)}\bigg(\frac{\mu_X^2}{\sigma_X^2}+\frac{\mu_Y^2}{\sigma_Y^2}-\frac{2\rho(x+\mu_X\mu_Y)}{\sigma_X\sigma_Y}\bigg)\bigg\}\nonumber\\
&\quad\times\sum_{n=0}^\infty\sum_{m=0}^{2n}\frac{x^{2n-m}|x|^{m-n}\sigma_X^{m-n-1}}{\pi(2n)!(1-\rho^2)^{2n+1/2}\sigma_Y^{m-n+1}}\binom{2n}{m}\bigg(\frac{\mu_X}{\sigma_X^2}-\frac{\rho \mu_Y}{\sigma_X\sigma_Y}\bigg)^m\nonumber\\
\label{pdf}&\quad\times\bigg(\frac{\mu_Y}{\sigma_Y^2}-\frac{\rho \mu_X}{\sigma_X\sigma_Y}\bigg)^{2n-m}K_{m-n}\bigg(\frac{|x|}{(1-\rho^2)\sigma_X\sigma_Y}\bigg), \quad x\in\mathbb{R},
\end{align}
where $K_\nu(x)=\int_0^\infty \mathrm{e}^{-x\cosh(t)}\cosh(\nu t)\,\mathrm{d}t$ is a modified Bessel function of the second kind \cite{olver}. 
If one of the means is equal to zero and $\rho=0$, then the density (\ref{pdf}) simplifies to a single infinite series (see \cite{cui,simon}). Without loss of generality, let $\mu_Y=0$. Then the PDF (\ref{pdf}) reduces to 
\begin{equation}\label{simple}
p_Z(x)=\frac{1}{\pi\sigma_X\sigma_Y}\exp\bigg(\!-\frac{\mu_X^2}{2\sigma_X^2}\bigg)\sum_{n=0}^\infty \frac{\mu_X^{2n}|x|^n}{(2n)!\sigma_X^{3n}\sigma_Y^n} K_n\bigg(\frac{|x|}{\sigma_X\sigma_Y}\bigg),\quad x\in\mathbb{R}.   
\end{equation}
In the case $\mu_X=\mu_Y=0$, 
the following exact formula for the PDF of $\overline{Z}_n$ was derived independently by \cite{man,np16,wb32}:
\begin{equation*}p_{\overline{Z}_n}(x)=\frac{2^{(1-n)/2}|x|^{(n-1)/2}}{s_n^{(n+1)/2}\sqrt{\pi(1-\rho^2)}\Gamma(n/2)}\exp\bigg(\frac{\rho  x}{s_n(1-\rho^2)} \bigg)K_{\frac{n-1}{2}}\bigg(\frac{ |x|}{s_n(1-\rho^2)}\bigg), \quad x\in\mathbb{R},
\end{equation*}
where $s_n=\sigma_X\sigma_Y/n$, and it was identified by \cite{gaunt prod} that $\overline{Z}_n$ is variance-gamma distributed. 
An exact formula for the PDF of the mean $\overline{Z}_n$ was recently obtained by \cite{gnp24} for the case of a general mean vector $(\mu_X,\mu_Y)\in\mathbb{R}^2$.

Over the years, other key distributional properties have been established for the distribution of $Z$, such as formulas for the characteristic function \cite{craig}, formulas for cumulants \cite{craig} and lower order moments \cite{h42}, location of the mode \cite{gz23} as well asymptotic approximations for the PDF, cumulative distribution function and quantile function \cite{gz23}. A review of the basic distributional theory for the distribution of $\overline{Z}_n$ in the case $\mu_X=\mu_Y=0$ is given in \cite{gaunt22}. 


In this paper, we establish another fundamental distributional property, by obtaining a Stein characterisation for the distribution $Z$, and more generally of $\overline{Z}_n$, in the general case $(\mu_X,\mu_Y)\in\mathbb{R}^2$, $(\sigma_X^2,\sigma_Y^2)\in(0,\infty)^2$, $-1<\rho<1$, $n\geq1$ (see Theorem \ref{cor2.2}). That is, we find a \emph{Stein operator} $A$ acting on a class of functions $\mathcal{F}$ such that $\mathbb{E}[Af(X)]=0$ for all $f\in\mathcal{F}$ if and only if $X=_d\overline{Z}_n$.  
Stein's classical characterisation of the normal distribution \cite{stein} states that $X\sim N(\mu,\sigma^2)$ if and only if
\begin{equation}\label{nsc}\mathbb{E}[\sigma^2f'(X)-(X-\mu)f(X)]=0
\end{equation}
for all absolutely continuous $f:\mathbb{R}\rightarrow\mathbb{R}$ such that $\mathbb{E}|f'(Y)|<\infty$ for $Y\sim N(\mu,\sigma^2)$. Since Stein's seminal work \cite{stein}, Stein characterisations have been established for many classical probability distributions (see \cite{agg1,gms,ley,mrs21} for an overview of this rather large literature), as well as more exotic distributions, such as stable distributions \cite{ah19,liu}, products of independent random variables \cite{gaunt pn,gms} and linear combinations of gamma random variables \cite{arras}, for which the PDF is difficult to write down in exact form. 
Our Stein characterisation generalises a number of Stein characterisations from this literature (see Remark \ref{genrem}), and has the interesting feature that the order of the characterising Stein operator decreases from fourth order (we refer to this as a fourth order Stein operator) in the general case $(\mu_X,\mu_Y)\in\mathbb{R}^2$ to second order in the case $\mu_X=\mu_Y=0$; see the Stein operator (\ref{zerozero}). There is in fact an intermediate level of complexity in that in the case $\mu_X/\sigma_X=\mu_Y/\sigma_Y$ there exists a third order Stein operator for $\overline{Z}_n$ (see Theorem \ref{cor2.4}). The orders of these operators are optimal amongst all Stein operators with linear coefficients; see Remark \ref{min}. Derivations of the Stein characterisations of Theorems \ref{cor2.2} and \ref{cor2.4} are given in Section \ref{sec3}.

Stein characterisations have classically been used as part of Stein's method to establish distributional approximations, finding applications in fields as diverse as random graph theory \cite{bhj92}, queuing theory \cite{bd16} and number theory \cite{h09}.
However, in recent years, there has been a growing trend in which Stein characterisations have found applications beyond proving quantitative limit theorems, with examples including new tests for goodness-of-fit \cite{steinstat}, new methods for parameter estimation \cite{ebner,fgs23}, relaxing the Gaussian assumption in shrinkage \cite{fathi} and derivations of distributional properties \cite{gms}.

In Section \ref{sec2.2}, 
we apply our Stein characterisations to obtain recurrence relations for the raw and central moments of $Z$ and $\overline{Z}_n$, which allow for efficient computation of higher order moments. From our recursions, we obtain, what are to the best of our knowledge, new formulas for the first four moments of $\overline{Z}_n$ from which we deduce new formulas for the skewness and kurtosis. 

\vspace{3mm}

\noindent \emph{Notation.} 
To simplify expressions, we define $r_X=\mu_X/\sigma_X$, $r_Y=\mu_Y/\sigma_Y$, $s_n=\sigma_X\sigma_Y/n$ and $s:=s_1=\sigma_X\sigma_Y$. 

\section{Main results}
\subsection{Stein characterisations}
In  the following Theorems \ref{cor2.2} and \ref{cor2.4}, we provide Stein characterisations of the mean $\overline{Z}_n=n^{-1}\sum_{i=1}^n Z_i$, where $Z_1,\ldots,Z_n$ are independent copies of $Z\sim PN(\mu_X,\mu_Y;\sigma_X^2,\sigma_Y^2;\rho)$. Note that setting $n=1$ yields Stein characterisations for $Z\sim PN(\mu_X,\mu_Y;\sigma_X^2,\sigma_Y^2;\rho)$.

\begin{theorem}\label{cor2.2} Let $\mu_X,\mu_Y\in\mathbb{R}$, $\sigma_X,\sigma_Y>0$ and $\rho\in(-1,1)$. Let $W$ be a real-valued random variable such that $\mathbb{E}|W|<\infty$. Define the operator $A_1$ by
\begin{align}
		A_1f(x)&=s_n^4(1-\rho^2)^2xf^{(4)}(x)+s_n^3(1-\rho^2)\big(ns_n(1-\rho^2)+4\rho  x\big)f^{(3)}(x)\nonumber\\
		&\quad+s_n^2\big(ns_n(\rho (r_X^2+r_Y^2)-(1+\rho^2) r_Xr_Y+3\rho(1-\rho^2))+(6\rho^2-2)x\big)f''(x)\nonumber\\
		&\quad+s_n\big(ns_n(2\rho r_Xr_Y-r_X^2-r_Y^2+3\rho^2-1)-4\rho  x\big)f'(x)\nonumber\\
		&\quad+\big(x-\mu_X\mu_Y-ns_n\rho\big)f(x).\label{steinoperator2}
	\end{align}
Then $W=_d\overline{Z}_n$ if and only if 
	\begin{equation}
		\mathbb{E}[A_1f(W)]=0\label{steinc}
	\end{equation}
for all $f:\mathbb{R}\rightarrow\mathbb{R}$ such that $f\in C^4(\mathbb{R})$ and $\mathbb{E}|f^{(m)}(\overline{Z}_n)|$, for $m=0,1,2,3$, and $\mathbb{E}|\overline{Z}_nf^{(m)}(\overline{Z}_n)|$, for $m=0,1,2,3,4$, are finite, where $f^{(0)}\equiv f$.
\end{theorem}

\begin{remark}
In the degenerate case $\rho\in\{-1,1\}$, the distribution of the product $Z=XY$ follows a non-central chi-square distribution. In this case, the Stein operator (\ref{steinoperator2}) reduces to    
\begin{align}
A_1f(x)&= s_n^2\big(ns_n(\rho (r_X^2+r_Y^2)-2 r_Xr_Y)+4x\big)f''(x)\nonumber\\
		&\quad+s_n\big(ns_n(2\rho r_Xr_Y-r_X^2-r_Y^2+2)-4\rho  x\big)f'(x)+\big(x-\mu_X\mu_Y-ns_n\rho\big)f(x),  \label{fd}
\end{align}
and the Stein characterisation reads:
$W=_d\overline{Z}_n$ if and only if 
	\begin{equation}
		\mathbb{E}[A_1f(W)]=0\label{steincaa}
	\end{equation}
for all $f:\mathbb{R}\rightarrow\mathbb{R}$ such that $f\in C^2((0,\infty))$ for $\rho=1$ (and $f\in C^2((-\infty,0))$ for $\rho=-1$) and $\mathbb{E}|f^{(m)}(\overline{Z}_n)|$, for $m=0,1,2$, and $\mathbb{E}|\overline{Z}_nf^{(m)}(\overline{Z}_n)|$, for $m=0,1,2$, are finite. 

We remark that a Stein operator for the non-central chi-square distribution, which is in agreement with the Stein operator (\ref{fd}) when $n=1$, has been obtained by \cite[Propositions 2.1 and 2.3]{g34}. The Stein operator (in the case $n=1$) was shown to be characterising in the sense of the Stein characterisation (\ref{steincaa}) in the first arXiv version of \cite{g34}, and the method we use to prove Theorem \ref{cor2.2} can also be used to prove this for the general $n\geq1$ case. Indeed, we stress that our proof of Theorem \ref{cor2.2} can be extended to deal with the degenerate case $\rho=\{-1,1\}$ with very minor modifications, and the reason we do not state Theorem \ref{cor2.2} for the full level of generality $\rho\in[-1,1]$, electing instead to state it in this remark, is because the class of functions under which the characterisation holds is different in the degenerate and non-degenerate cases.
\end{remark}

In the case $\mu_X/\sigma_X=\mu_Y/\sigma_Y$, we can reduce the order of the Stein operator for $\overline{Z}_n$ from fourth order to third order.

\begin{theorem}\label{cor2.4} Suppose that $\mu_X/\sigma_X=\mu_Y/\sigma_Y$. Let $W$ be a real-valued random variable such that $\mathbb{E}|W|<\infty$.  Define the operator $A_2$ by
\begin{align}\label{steinoperator4}
		A_{2}f(x)&=s_n^3(1-\rho^2)(1+\rho)xf^{(3)}(x)+s_n^2(1+\rho)\big(ns_n(1-\rho^2)+(3\rho-1)x\big)f''(x)\nonumber\\
		&\quad+s_n\big(ns_n(2\rho^2+\rho-1-(1-\rho)r_Xr_Y)-(3\rho+1)x\big)f'(x)\nonumber\\
  &\quad+\big(x-ns_n\rho-\mu_X\mu_Y\big)f(x).
	\end{align}
Then $W=_d\overline{Z}_n$ if and only if 
	\begin{equation}
		\mathbb{E}[A_{2}f(W)]=0\label{charac4}
	\end{equation}
for all $f:\mathbb{R}\rightarrow\mathbb{R}$ such that $f\in C^3(\mathbb{R})$ and $\mathbb{E}|f^{(m)}(\overline{Z}_n)|$, for $m=0,1,2$, and $\mathbb{E}|\overline{Z}_nf^{(m)}(\overline{Z}_n)|$, for $m=0,1,2,3$, are finite.
\end{theorem}

\begin{remark}Suppose that $\mu_X/\sigma_X=\mu_Y/\sigma_Y$. We note that, in this case, making the substitution 
$g(x)=(1-\rho)s_nf'(x)+f(x)$ in the Stein operator (\ref{steinoperator2}) yields  
\begin{align*}
		A_{1}f(x)&=s_n^3(1-\rho^2)(1+\rho)xg^{(3)}(x)+s_n^2(1+\rho)\big(ns_n(1-\rho^2)+(3\rho-1)x\big)g''(x)\\
		&\quad+s_n\big(ns_n(2\rho^2+\rho-1-(1-\rho)r_Xr_Y)-(3\rho+1)x\big)g'(x)\\
  &\quad+\big(x-ns_n\rho-\mu_X\mu_Y\big)g(x),
	\end{align*}
which we recognise as the Stein operator (\ref{steinoperator4}).
\end{remark}

\begin{remark}\label{genrem} The Stein operator (\ref{steinoperator2}) reduces to a number of other Stein operators from the literature as special cases. Taking $\mu_X=\mu_Y=0$ in (\ref{steinoperator2}), we obtain
\begin{align}
		A_3f(x)&=s_n^4(1-\rho^2)^2xf^{(4)}(x)+(1-\rho^2)s_n^3(ns_n(1-\rho^2)+4\rho x)f^{(3)}(x)+s_n^2((6\rho^2-2)x\nonumber\\
		&\quad+3ns_n\rho(1-\rho^2))f''(x)+s_n(ns_n(3\rho^2-1)-4\rho x)f'(x)+(x-ns_n\rho)f(x).\label{near}
	\end{align}
On making the substitution $g(x)=(1-\rho^2)s_n^2f''(x)+2\rho s_nf'(x)-f(x)$ in (\ref{near}), we obtain the Stein operator
	\begin{equation}\label{zerozero}
		A_4g(x)=s_n^2(1-\rho^2)xg''(x)+s_n(ns_n(1-\rho^2)+2\rho x)g'(x)+(ns_n\rho-x)g(x),
	\end{equation}
	which is given in \cite{gaunt22} and is a special case of the variance-gamma Stein operator of \cite{gaunt vg}; recall that $\overline{Z}_n$ is variance-gamma distributed when $\mu_X=\mu_Y=0$. Further setting $n=1$ and $\rho=0$ in (\ref{zerozero}) yields the Stein operator
 \begin{equation*}A_5g(x)=s^2xg''(x)+s^2g'(x)-xg(x),
 \end{equation*}
 which is a special case of the Stein operator of \cite{gaunt pn} for the product of $k\geq2$ independent zero mean normal random variables.
	
	 Now, for the sake of simpler expressions, we let $\sigma_X=\sigma_Y=1$. Setting $n=1$ and $\rho=0$ in (\ref{steinoperator2}) yields the Stein operator 
	\begin{align*}
		A_6f(x)=xf^{(4)}(x)+f^{(3)}(x)-(\mu_X\mu_Y+2x)f''(x)-(\mu^2_X+\mu^2_Y+1)f'(x)+(x-\mu_X\mu_Y)f(x),
	\end{align*}
which we recognise as the Stein operator of \cite{gms20} for the product of two independent normal random variables with general means $(\mu_X,\mu_Y)\in\mathbb{R}^2$. Now suppose $\mu_X=\mu_Y=\mu$ (and again we set $\sigma_X=\sigma_Y=1$). Then, setting $n=1$ and $\rho=0$ in (\ref{steinoperator4}) yields the Stein operator
\begin{align*}
	A_{7}f(x)=xf^{(3)}(x)+(1-x)f''(x)-(x+1+\mu^2)f'(x)+(x-\mu^2)f(x),
\end{align*}
which we identify as the Stein operator of \cite{gms20} for the product of two independent normal random variables with equal but possibly non-zero means.
\end{remark}

\begin{remark}\label{min}
It was shown in \cite{gms20} that there does not exist a third order Stein operator with linear coefficients for the product of two independent normal random variables with general means $(\mu_X,\mu_Y)\in\mathbb{R}^2$. It was also shown in \cite{gms20} that there does not exist a second order Stein operator with linear coefficients for the product of two independent non-centred normal random variables such that $\mu_X=\mu_Y$ (when $\sigma_X=\sigma_Y=1$). (Note that in this paper we have generalised their condition to $\mu_X/\sigma_X=\mu_Y/\sigma_Y$.) Thus, the order of the Stein operators (\ref{steinoperator2}) and (\ref{steinoperator4}) are optimal amongst all Stein operators with linear coefficients.
  
  In the case $\rho=0$ and $\mu_Y=0$, the formula (\ref{pdf}) for the PDF of $Z$ simplifies from a double infinite series of modified Bessel functions of the second kind to the single infinite series (\ref{simple}). It is therefore natural to ask whether it is possible to obtain a Stein operator in this case with linear coefficients and order strictly less than four (the analysis of \cite{gms20} outlined above did not rule out this possibility). However, we found that there is no reduction in the order of the Stein operator in this case. To show this, we will use the brute force approach of \cite{gms20}, to show that there is no third order Stein operator with linear coefficients for the product $Z=XY$, where $X\sim N(1,1)$ and $Y\sim N(0,1)$ are independent.
  Suppose for contradiction that there exists a third order Stein operator for $Z$, then it would be of the form $A_Zf(x)=\sum_{j=0}^{3}(a_{0,j}+a_{1,j}x)f^{(j)}(x)$, where $f^{(0)}\equiv f$. Now, if $A_Z$ was a Stein operator for $Z$, we would have $\mathbb{E}[A_Zf(Z)]=0$ for all $f$ in some class of functions that contains the monomials $\{x^k:k\geq1\}$. Taking $f(x)=x^k$, $k=0,1,\ldots,7$, we obtain eight equations for eight unknowns. Letting $\mu'_k$ denote $\mathbb{E}[Z^k]$, we get $\mu_1'=0$, $\mu_2'=2$, $\mu_3'=0$, $\mu_4'=30$, $\mu_5'=0$, $\mu_6'=1140$, $\mu_7'=0$, $\mu_8'=80220$. This yields the system of equations
  \begin{align*}
  	a_{0,0}&=0\\
  	2a_{1,0}+a_{0,1}&=0\\
  	2a_{0,0}+4a_{1,1}+2a_{0,2}&=0\\
  	30a_{1,0}+6a_{0,1}+12a_{1,2}+6a_{0,3}&=0\\
  	30a_{0,0}+120a_{1,1}+24a_{0,2}+48a_{1,3}&=0\\
  	1140a_{1,0}+150a_{0,1}+600a_{1,2}+120a_{0,3}&=0\\
  	1140a_{0,0}+6840a_{1,1}+900a_{0,2}+3600a_{1,3}&=0\\
  	80220a_{1,0}+7980a_{0,1}+47880a_{1,2}+6300a_{0,3}&=0.
  	\end{align*}
We used \emph{Mathematica} to compute that the determinant of the matrix corresponding to this system of linear equations is $125\,411\,328\,000\not=0$. Therefore, there is a unique solution which is $a_{0,0}=a_{1,0}=\cdots=a_{1,3}=0$. 
Therefore, in general, there does not exist a third order Stein operator with linear coefficients for a product of two independent normal random variables, where one of the means is zero and the other is non-zero.    
\end{remark}

\begin{remark}
In the zero mean case $\mu_X=\mu_Y=0$, the sample mean $\overline{Z}_n$ is variance-gamma distributed, and one can then write down a corresponding Stein equation to the Stein operator (\ref{zerozero}), which is a special case of the variance-gamma Stein equation of \cite{gaunt vg}. The variance-gamma Stein equation was solved by \cite{gaunt vg} and bounds for the derivatives of the solution of the Stein equation are given in \cite{dgv1,vgd2,vgd3}. However, in the non-zero mean case we have been unable to even write down a solution to a Stein equation corresponding to the Stein operator (\ref{steinoperator2}). We consider it to be an interesting but very difficult problem to obtain bounds on the derivatives of the solution of the Stein equation corresponding to the Stein operator (\ref{steinoperator2}).
\end{remark}

\subsection{Applications of the Stein characterisations}\label{sec2.2}

\subsubsection{Moments}\label{2.1.1}

Substituting $f_1(x)=x^k$ and $f_2(x)=(x-\mathbb{E}[\overline{Z}_n])^k=(x-(\mu_X\mu_Y+ns_n\rho))^k$ into the characterising equation (\ref{steinc}) yields the following recursions for the $k$-th raw moment $\mu'_k=\mathbb{E}[\overline{Z}_n^k]$ and the $k$-th central moment $\mu_k=\mathbb{E}[(\overline{Z}_n-\mathbb{E}[\overline{Z}_n])^k]$:
\begin{align}
	\mu_{k+1}'&=\big(\mu_X\mu_Y+s_n\rho(4k+n)\big)\mu_k'-s_n^2k\big(n(2\rho r_Xr_Y-r_X^2-r_Y^2+3\rho^2-1)\nonumber\\
	&\quad+(k-1)(6\rho^2-2)\big)\mu_{k-1}'-s_n^3k(k-1)\big(n(\rho r_Y^2-\rho^2r_Xr_Y-r_Xr_Y+\rho r_X^2+3\rho(1-\rho^2))\nonumber\\
	&\quad+4(k-2)\rho(1-\rho^2)\big)\mu_{k-2}'-s_n^4k(k-1)(k-2)(1-\rho^2)^2(n+k-3)\mu_{k-3}',\quad\text{$k\geq0$},\label{rec1}\\
	\mu_{k+1}&=4\rho s_nk\mu_k-s_nk\big(s_n(6\rho^2-2)(k-1)+ns_n(2\rho r_Xr_Y-r_X^2-r_Y^2+3\rho^2-1)\nonumber\\
&\quad-4\rho(\mu_X\mu_Y+ns_n\rho)\big)\mu_{k-1}-s_n^2k(k-1)\big((6\rho^2-2)(\mu_X\mu_Y+ns_n\rho)\nonumber\\
	&\quad+ns_n(\rho r_Y^2-\rho^2 r_Xr_Y-r_Xr_Y+\rho r_X^2+3\rho(1-\rho^2))+4s_n\rho(1-\rho^2)(k-2)\big)\mu_{k-2}\nonumber\\
	&\quad-s_n^3(1-\rho^2)k(k-1)(k-2)\big(4\rho(\mu_X\mu_Y+ns_n\rho)+s_n(1-\rho^2)(n+k-3)\big)\mu_{k-3}\nonumber\\
	&\quad-s_n^4(1-\rho^2)^2k(k-1)(k-2)(k-3)(\mu_X\mu_Y+ns_n\rho)\mu_{k-4}, \quad\text{$k\geq0$}.\label{rec2}
\end{align}
These recurrences allow raw and central moments of general order to be efficiently computed, using just the basic starting values $\mu_0'=1$, $\mu_0=1$ and $\mu_1=0$. In the case  $\mu_X/\sigma_X=\mu_Y/\sigma_Y$, we can obtain lower order recursions: substituting $f_1(x)=x^k$ and $f_2(x)=(x-(\mu_X\mu_Y+ns_n\rho))^k$ into (\ref{charac4}) yields the recurrence relations
\begin{align*}
	\mu_{k+1}'&=\big(\mu_X\mu_Y+s_n\big(\rho n +(3\rho+1)k\big)\big)\mu_k'\\
	&\quad-s_nk\big(\mu_X\mu_Y(\rho-1)+s_n\big(n(2\rho^2+\rho-1)+(1+\rho)(3\rho-1)(k-1)\big)\big)\mu_{k-1}'\\
	&\quad-s_n^3(1+\rho)(1-\rho^2)k(k-1)(k-2+n)\mu_{k-2}',\quad\text{$k\geq0$},\\
	\mu_{k+1}&=ks_n(3\rho+1)\mu_k-ks_n\big((1+\rho)(3\rho-1)(k-1)s_n+\mu_X\mu_Y(\rho-1)+ns_n(2\rho^2+\rho-1)\\
	&\quad-(3\rho+1)(\mu_X\mu_Y+ns_n\rho)\big)\mu_{k-1}-s_n^2(1+\rho)k(k-1)\big((\mu_X\mu_Y+ns_n\rho)(3\rho-1)\\
	&\quad+n(1-\rho^2)s_n+(1-\rho^2)s_n(k-2)\big)\mu_{k-2}\\
	&\quad-s_n^3k(k-1)(k-2)(1-\rho^2)(1+\rho)(\mu_X\mu_Y+ns_n\rho)\mu_{k-3},\quad\text{$k\geq0$}.
\end{align*}
In the case $\mu_X=\mu_Y=0$, the recurrence relations for the raw and central moments simplify further, since in this case the Stein operator (\ref{zerozero}) is a second order differential operator; see \cite{gaunt22} for these simpler recursions.

Let us now consider the general case $(\mu_X,\mu_Y)\in\mathbb{R}^2$, $\sigma_X,\sigma_Y>0$ and $-1<\rho<1$. From the recursions (\ref{rec1}) and (\ref{rec2}), we readily obtain that the first four raw moments are given by 
\begin{align*}
	\mu_1'&=\mu_X\mu_Y+\rho\sigma_X\sigma_Y,\\
	\mu_2'&=\frac{\sigma_X^2\sigma_Y^2}{n}\Big\{nr_X^2r_Y^2+r_X^2+r_Y^2+2\rho(n+1)r_Xr_Y+\rho^2(n+1)+1\Big\},\\
	\mu_3'&=\frac{\sigma_X^3\sigma_Y^3}{n^2}\Big\{n^2r_X^3r_Y^3+3nr_Xr_Y(r_X^2+r_Y^2)+3\rho n(n+2)r_X^2r_Y^2+3\rho(n+2)(r_X^2+r_Y^2)\\
 &\quad+3(n+2)\big(\rho^2(n+1)+1\big)r_Xr_Y
	+\rho(n+2)\big(\rho^2(n+1)+3\big)\Big\},\\
	\mu_4'&=\frac{\sigma_X^4\sigma_Y^4}{n^3}\Big\{n^3r_X^4r_Y^4+4\rho n^2(n+3)r_X^3r_Y^3+6n^2r_X^2r_Y^2(r_X^2+r_Y^2)+3n(r_X^4+r_Y^4)\\
	&\quad+12\rho n(n+3)r_Xr_Y(r_X^2+r_Y^2)+6n\big(\rho^2(n+2)(n+3)+(n+5)\big)r_X^2r_Y^2\\
 &\quad+6(n+2)\big(\rho^2(n+3)+1\big)(r_X^2+r_Y^2)+4\rho(n+2)(n+3)\big(\rho^2(n+1)+3\big)r_Xr_Y\\
	&\quad+\rho^4(n+1)(n+2)(n+3)+6\rho^2(n+2)(n+3)+3(n+2)\Big\},
\end{align*}
and that the first four central moments are given by 
\begin{align}
	\mu_1&=0,\nonumber\\
	\mu_2&=\frac{\sigma_X^2\sigma_Y^2}{n}\Big\{r_X^2+r_Y^2+2\rho r_Xr_Y+\rho^2+1\Big\},\label{2222}\\
	\mu_3&=\frac{2\sigma_X^3\sigma_Y^3}{n^2}\Big\{3\rho(r_X^2+r_Y^2)+3(\rho^2+1)r_Xr_Y+\rho(\rho^2+3)\Big\},\nonumber\\
	\mu_4&=\frac{3\sigma_X^4\sigma_Y^4}{n^3}\Big\{n(r_X^4+r_Y^4)+4\rho n r_Xr_Y(r_X^2+r_Y^2)+2n(2\rho^2+1)r_X^2r_Y^2\nonumber\\
	&\quad+2\big(\rho^2(n+6)+(n+2)\big)(r_X^2+r_Y^2)+4\rho\big(\rho^2(n+2)+(n+6)\big)r_Xr_Y\nonumber\\
	&\quad+\rho^4(n+2)+2\rho^2(n+6)+(n+2)\Big\}.\label{4444}
\end{align}
The variance is given by 
\[\mathrm{Var}(\overline{Z}_n)=\mu_2=\frac{\sigma_X^2\sigma_Y^2}{n}\Big\{r_X^2+r_Y^2+2\rho r_Xr_Y+\rho^2+1\Big\},\]
and the skewness $\gamma_1=\mu_3/\mu_2^{3/2}$ and kurtosis $\beta_2=\mu_4/\mu_2^2$ can also be computed using the above formulas.

To the best of our knowledge, the above formulas for the moments of $\overline{Z}_n$ are new. Formulas for the first four raw moments of the product $Z$ are given in \cite{h42} and formulas for the variance, skewness and kurtosis are given in \cite{skewkurt}. However, we noticed that the formula of \cite[Proposition 2.2]{skewkurt} for the kurtosis of $Z$ is incorrect; when our formulas (\ref{2222}) and (\ref{4444}) for the second and fourth central moments of $\overline{Z}_n$ are specialised to the case $n=1$ and the kurtosis is calculated using $\beta_2=\mu_4/\mu_2^2$, the resulting formula corrects their formula for the kurtosis of $Z$, which we report here:
\begin{align*}
\mathrm{Kurt}[Z]&=\frac{3\big(r_X^4+r_Y^4+4\rho  r_Xr_Y(r_X^2+r_Y^2)+2(2\rho^2+1)r_X^2r_Y^2\big)}{(r_X^2+r_Y^2+2\rho r_Xr_Y+\rho^2+1)^2}\\
	&\quad+\frac{3\big(2(7\rho^2+3)(r_X^2+r_Y^2)+4\rho(3\rho^2+7)r_Xr_Y+3\rho^4+14\rho^2+3\big)}{(r_X^2+r_Y^2+2\rho r_Xr_Y+\rho^2+1)^2}.
\end{align*}

\subsubsection{Probability density function}\label{2.1.3}

An explicit formula for the PDF $p_n:=p_{\overline{Z}_n}$ of $\overline{Z}_n$ was recently obtained by \cite{gnp24}. Here, we apply our Stein characterisation (\ref{steinc}) to derive a differential equation satisfied by the PDF $p_n$, which provides some insight into the complexity of this formula. To do this, we recall a duality result of \cite[Remark 2.7]{gms} (see also Section 4 of that paper for further details). Define the operators $M(f):x\mapsto xf(x)$ and $D(f):x\mapsto f'(x)$. If $V$ has a smooth density $p$, which solves the ODE $B p = 0$ with $B = \sum_{i,j} b_{ij} M^j D^i$, where the $b_{ij}$ are real-valued constants, then a Stein operator for $V$ is given by $A_V = \sum_{i,j} (-1)^i b_{ij} D^i M^j$, and similarly given a Stein operator for $V$ one can write down an ODE that $p$ satisfies. In this manner, we can obtain the following ODE satisfied by the density $p_n$ of $\overline{Z}_n$ (where we set $\sigma_X=\sigma_Y=1$ for sake of a simpler expression): 
\begin{align}
	&(1-\rho^2)^2xp_{n}^{(4)}(x)+(1-\rho^2)\big((1-\rho^2)(4-n)-4\rho n x\big)p^{(3)}_{n}(x)\nonumber\\
	&\quad+n\big((6\rho^2-2)nx-12\rho(1-\rho^2)+n(3\rho(1-\rho^2)+\rho \mu_Y^2-\rho^2 \mu_X\mu_Y-\mu_X\mu_Y+\rho \mu_X^2)\big)p''_{n}(x)\nonumber\\
	&\quad+n^2\big(2(6\rho^2-2)-n(2\rho \mu_X\mu_Y-\mu_X^2-\mu_Y^2+3\rho^2-1)+4\rho nx\big)p'_{n}(x)\nonumber\\
	&\quad+n^3\big(4\rho+n(x-\mu_X\mu_Y-\rho)\big)p_{n}(x)=0. \label{ODE1}
	\end{align}

From Remark \ref{min} and the duality result from \cite{gms} outlined above, we know there does not exist an ODE satisfied by $p_{n}$ with linear coefficients with a lower degree than (\ref{ODE1}). 
We were unable to transform (\ref{ODE1})  into a known class such as the Meijer-$G$ function differential equation, therefore it seems plausible that the formula (\ref{pdf}) for the PDF of $Z$ cannot be simplified further.
We also note that there is not a severe increase in complexity in the ODE (\ref{ODE1}) satisfied by $\overline{Z}_n$ from the $n=1$ case to the general $n\geq1$ case. 

\section{Proofs}\label{sec3}

For sake of simpler expressions, we first derive our Stein characterisations for the case $\sigma_X=\sigma_Y=1$, and then generalise to the general case $\sigma_X,\sigma_Y>0$ by using the distributional relation $Z=XY=_d \sigma_X\sigma_Y UV$, where $(U,V)$ is a bivariate normal random vector with mean vector $(\mu_X/\sigma_X,\mu_Y/\sigma_Y)$, variances $(1,1)$ and correlation coefficient $\rho$.

\vspace{3mm}

\noindent{\emph{Proof of Theorem \ref{cor2.2}.}} \emph{Necessity}. Suppose that $\sigma_X=\sigma_Y=1$. 
We start by
defining the random variable $V$ by $V=(Y-\rho X)/\sqrt{1-\rho^2}$, so that $X$ and $V$ are independent standard normal random variables. We can write $Z$ as 
\begin{align*}
	Z&=(X+\mu_X)(\rho X+\sqrt{1-\rho^2}V+\mu_Y).
\end{align*}
Therefore, $Z|X\sim N((X+\mu_X)(\rho X+\mu_Y), (1-\rho^2)(X+\mu_X)^2)$, and we obtain from (\ref{nsc}) that 
\begin{align}
	\mathbb{E}[Zf(Z)]&=\mathbb{E}\big[\mathbb{E}[Zf(Z)|X]\big]\nonumber\\
	&=\mathbb{E}\big[\mathbb{E}[(1-\rho^2)(X+\mu_X)^2f'(Z)+(X+\mu_X)(\rho X+\mu_Y)f(Z)|X]\big]\nonumber\\
	&=\mathbb{E}[(1-\rho^2)(X+\mu_X)^2f'(Z)]+\mathbb{E}[\rho X(X+\mu_X)f(Z)]+\mathbb{E}[\mu_Y(X+\mu_X)f(Z)].\label{star}
\end{align} 
Now, let $z=(x+\mu_X)(\rho x+\sqrt{1-\rho^2}v+\mu_Y)$ and note that
\begin{align}
	\frac{\partial z}{\partial x}
	&=2\rho(x+\mu_X)+(\mu_Y-\rho\mu_X)+\sqrt{1-\rho^2}v,\label{partial3}\\
	(x+\mu_X)\frac{\partial z}{\partial x}&=z+\rho(x+\mu_X)^2\label{partial1},\\
		\frac{\partial z}{\partial v}&=\sqrt{1-\rho^2}(x+\mu_X).\label{partial2}
\end{align}
We will also use repeatedly that, for $X\sim N(0,1)$,
\begin{equation}
	\mathbb{E}[g'(X)-Xg(X)]=0\label{snsc},
\end{equation}
for all differentiable $g:\mathbb{R}\rightarrow\mathbb{R}$ such that $\mathbb{E}
|g'(Y)|$ is finite for $Y\sim N(0,1)$.

Now, by conditioning on $V$, we obtain
\begin{align}
	\mathbb{E}[Zf(Z)]&=\mathbb{E}\big[\mathbb{E}[\rho X\cdot(X+\mu_X)f(Z)|V]\big]+\mathbb{E}[\mu_Y(X+\mu_X)f(Z)]\nonumber\\
	&\quad+\mathbb{E}[(1-\rho^2)(X+\mu_X)^2f'(Z)]\nonumber\\
	&=\mathbb{E}[\rho f(Z)]+\mathbb{E}[\rho Zf'(Z)]+\mathbb{E}[\rho^2(X+\mu_X)^2f'(Z)]\nonumber\\
	&\quad+\mathbb{E}[\mu_Y(X+\mu_X)f(Z)]+\mathbb{E}[(1-\rho^2)(X+\mu_X)^2f'(Z)]\nonumber,
\end{align}
which follows by using (\ref{snsc}) with $g(x)=\rho(x+\mu_X)f(z)$, together with (\ref{partial1}).

By conditioning on $V$ again, we have
\begin{align*}
	\mathbb{E}[Zf(Z)]&=\mathbb{E}[\rho f(Z)]+\mathbb{E}[\rho Zf'(Z)]+\mathbb{E}[\rho^2(X+\mu_X)^2f'(Z)]+\mathbb{E}[\mu_Y(X+\mu_X)f(Z)]\\
	&\quad+\mathbb{E}\big[\mathbb{E}[(1-\rho^2)X\cdot(X+\mu_X)f'(Z)|V]\big]+\mathbb{E}[(1-\rho^2)\mu_X(X+\mu_X)f'(Z)]\\
	&=\mathbb{E}[\rho f(Z)]+\mathbb{E}[\rho Zf'(Z)]+\mathbb{E}[\rho^2(X+\mu_X)^2f'(Z)]\\
	&\quad+\mathbb{E}[\mu_Y(X+\mu_X)f(Z)]+\mathbb{E}[(1-\rho^2)f'(Z)]+\mathbb{E}[(1-\rho^2)Zf''(Z)]\\
	&\quad+\mathbb{E}[(1-\rho^2)\rho(X+\mu_X)^2f''(Z)]+\mathbb{E}[(1-\rho^2)\mu_X(X+\mu_X)f'(Z)],
\end{align*}
which follows on using (\ref{snsc}) with $g(x)=(1-\rho^2)(x+\mu_X)f'(z)$, together with (\ref{partial1}).
This implies that 
\begin{align}
	&\mathbb{E}[(\rho-Z)f(Z)]+\mathbb{E}[\rho Zf'(Z)]+\mathbb{E}[\rho^2(X+\mu_X)^2f'(Z)]\nonumber\\
	&+\mathbb{E}[\mu_Y(X+\mu_X)f(Z)]+\mathbb{E}[(1-\rho^2)f'(Z)]+\mathbb{E}[(1-\rho^2)Zf''(Z)]\nonumber\\
	&+\mathbb{E}[(1-\rho^2)\rho(X+\mu_X)^2f''(Z)]+\mathbb{E}[(1-\rho^2)\mu_X(X+\mu_X)f'(Z)]=0.\label{eqn1}
\end{align}
If we now consider the equation (\ref{star}) multiplied by $\rho$ with $f$ replaced by $f'$, we have the equality 
\begin{align*}
	\mathbb{E}[\rho Zf'(Z)]&=\mathbb{E}[\rho^2X(X+\mu_X)f'(Z)]+\mathbb{E}[\rho\mu_Y(X+\mu_X)f'(Z)]\\
 &\quad+\mathbb{E}[(1-\rho^2)\rho(X+\mu_X)^2f''(Z)]\\
	&=\mathbb{E}[\rho^2(X+\mu_X)^2f'(Z)]+\mathbb{E}[\rho\mu_Y(X+\mu_X)f'(Z)]\\
	&\quad+\mathbb{E}[(1-\rho^2)\rho(X+\mu_X)^2f''(Z)]-\mathbb{E}[\rho^2\mu_X(X+\mu_X)f'(Z)].
\end{align*}
By rearranging it follows that
\begin{align}
	&\mathbb{E}[\rho^2(X+\mu_X)^2f'(Z)]+\mathbb{E}[(1-\rho^2)\rho(X+\mu_X)^2f''(Z)]\nonumber\\
	&=\mathbb{E}[\rho Zf'(Z)]+\mathbb{E}[\rho^2\mu_X(X+\mu_X)f'(Z)]-\mathbb{E}[\rho\mu_Y(X+\mu_X)f'(Z)].\label{eqn2}
\end{align}
Now, substituting (\ref{eqn2}) into (\ref{eqn1}) gives that
\begin{align}
	&\mathbb{E}[(1-\rho^2)Zf''(Z)]+\mathbb{E}[(1-\rho^2+2\rho Z)f'(Z)]+\mathbb{E}[(\rho-Z)f(Z)]\nonumber\\
	&+\mathbb{E}[\mu_Y(X+\mu_X)f(Z)]+\mathbb{E}[\mu_X(X+\mu_X)f'(Z)]-\mathbb{E}[\rho\mu_Y(X+\mu_X)f'(Z)]=0\label{F}.
\end{align}
Adding together (\ref{F}) multiplied by $(1-\rho^2)$ with $f$ replaced by $f''$ and (\ref{F}) multiplied by 2$\rho$ with $f$ replaced by $f'$, and subtracting (\ref{F}), we have 
\begin{align}
	&\mathbb{E}[(1-\rho^2)^2Zf^{(4)}(Z)]+\mathbb{E}[(1-\rho^2)(1-\rho^2+2\rho Z)f^{(3)}(Z)]+\mathbb{E}[(1-\rho^2)(\rho-Z)f''(Z)]\nonumber\\
 &+\mathbb{E}[2\rho(1-\rho^2)Zf^{(3)}(Z)]+\mathbb{E}[2\rho(1-\rho^2+2\rho Z)f''(Z)]+\mathbb{E}[2\rho(\rho-Z)f'(Z)]\nonumber\\
	&-\mathbb{E}[(1-\rho^2)Zf''(Z)]-\mathbb{E}[(1-\rho^2+2\rho Z)f'(Z)]\nonumber\\
	&-\mathbb{E}[(\rho-Z)f(Z)]+\mathbb{E}[(1-\rho^2)\mu_Y(X+\mu_X)f''(Z)]\nonumber\\
	&+\mathbb{E}[(1-\rho^2)(\mu_X-\rho\mu_Y)(X+\mu_X)f^{(3)}(Z)]+\mathbb{E}[2\rho\mu_Y(X+\mu_X)f'(Z)]\nonumber\\
	&+\mathbb{E}[2\rho(\mu_X-\rho\mu_Y)(X+\mu_X)f''(Z)]+\mathbb{E}[(X+\mu_X)(\rho\mu_Y-\mu_X)f'(Z)]\nonumber\\
	&-\mathbb{E}[\mu_Y(X+\mu_X)f(Z)]=0\label{eqn4}.
\end{align}
Now consider the final two terms from (\ref{eqn4}). By conditioning on $X$ and $V$ and using both (\ref{partial3}) and (\ref{partial2}), we obtain
\begin{align}
	&\mathbb{E}[(X+\mu_X)(\rho\mu_Y-\mu_X)f'(Z)]-\mathbb{E}[\mu_Y(X+\mu_X)f(Z)]\nonumber\\
	&=\mathbb{E}[(X+\mu_X)(\rho\mu_Y-\mu_X)f'(Z)]-\mathbb{E}[\mu_X\mu_Yf(Z)]-\mathbb{E}\big[\mathbb{E}[X\cdot\mu_Yf(Z)|V]\big]\nonumber\\
	&=\mathbb{E}[(X(\rho\mu_Y-\mu_X)f'(Z)]+\mathbb{E}[\mu_X(\rho\mu_Y-\mu_X)f'(Z)]-\mathbb{E}[\mu_X\mu_Yf(Z)]\nonumber\\
	&\quad-\mathbb{E}[2\rho\mu_Y(X+\mu_X)f'(Z)]-\mathbb{E}[\mu_Y(\mu_Y-\rho\mu_X)f'(Z)]-\mathbb{E}\big[\mathbb{E}[V\cdot\mu_Y\sqrt{1-\rho^2} f'(Z)|X]\big]\nonumber\\
	&=\mathbb{E}\big[\mathbb{E}[X\cdot(\rho\mu_Y-\mu_X)f'(Z)|V]\big]-\mathbb{E}[(\mu_X^2+\mu_Y^2)f'(Z)]-\mathbb{E}[\mu_X\mu_Yf(Z)]\nonumber\\
	&\quad-\mathbb{E}[2\rho\mu_Y Xf'(Z)]-\mathbb{E}[\mu_Y(1-\rho^2)(X+\mu_X)f''(Z)]\nonumber\\
	&=\mathbb{E}[2\rho(\rho\mu_Y-\mu_X)(X+\mu_X)f''(Z)]+\mathbb{E}[(\rho\mu_Y-\mu_X)(\mu_Y-\rho\mu_X)f''(Z)]\nonumber\\
	&\quad+\mathbb{E}\big[\mathbb{E}[V\cdot(\rho\mu_Y-\mu_X)\sqrt{1-\rho^2}f''(Z)|X]\big]-\mathbb{E}[(\mu^2_X+\mu_Y^2)f'(Z)]-\mathbb{E}[\mu_X\mu_Yf(Z)]\nonumber\\
	&\quad-\mathbb{E}[2\rho\mu_YXf'(Z)]-\mathbb{E}[\mu_Y(1-\rho^2)(X+\mu_X)f''(Z)]\nonumber\\
	&=\mathbb{E}[2\rho(\rho\mu_Y-\mu_X)(X+\mu_X)f''(Z)]+\mathbb{E}[(\rho\mu_Y-\mu_X)(\mu_Y-\rho\mu_X)f''(Z)]\nonumber\\
	&\quad+\mathbb{E}[(\rho\mu_Y-\mu_X)(1-\rho^2)(X+\mu_X)f^{(3)}(Z)]-\mathbb{E}[(\mu_X^2+\mu_Y^2)f'(Z)]-\mathbb{E}[\mu_X\mu_Yf(Z)]\nonumber\\
	&\quad-\mathbb{E}[2\rho \mu_YXf'(Z)]-\mathbb{E}[\mu_Y(1-\rho^2)(X+\mu_X)f''(Z)],\label{eqn5}
\end{align}
where the second, third, fourth and fifth equalities follow by using (\ref{snsc}) with $g(x)=\mu_Yf(z)$, $g(v)=\mu_Y\sqrt{1-\rho^2}f'(z)$, $g(x)=(\rho\mu_Y-\mu_X)f'(z)$ and $g(v)=(\rho\mu_Y-\mu_X)\sqrt{1-\rho^2}f''(z)$, respectively.
Finally, substituting (\ref{eqn5}) into (\ref{eqn4}) and rearranging, we obtain that, for $Z\sim PN(\mu_X,\mu_Y;1,1;\rho)$,
\begin{align}
	&\mathbb{E}\big[(1-\rho^2)^2Zf^{(4)}(Z)+(1-\rho^2)\big((1-\rho^2)+4\rho Z\big)f^{(3)}(Z)\nonumber\\
	&+\big((\rho(\mu_X^2+\mu^2_Y)-(1+\rho^2)\mu_X\mu_Y+3\rho(1-\rho^2))+(6\rho^2-2)Z\big)f''(Z)\nonumber\\
	&+\big((2\rho\mu_X\mu_Y-\mu^2_X-\mu_Y^2+3\rho^2-1)-4\rho Z\big)f'(Z)+(Z-\mu_X\mu_Y-\rho)f(Z)\big]=0,\nonumber
\end{align}
which we recognise as the characterising equation (\ref{steinc}) in the case $\sigma_X=\sigma_Y=1$ and $n=1$. 

We now extend to $n\geq1$. Let $S=n\overline{Z}_n=\sum_{i=1}^nZ_i$. Then, by conditioning,
\begin{align}&\mathbb{E}[ (S-n\mu_X\mu_Y-n\rho)f(S)]=\sum_{i=1}^n\mathbb{E}\big[\mathbb{E}\big[(Z_i-\mu_X\mu_Y-\rho)f(S)\,\big|\,\{Z_j\}_{j\not=i}\big]\big]\nonumber\\
&=-\sum_{i=1}^n\mathbb{E}\big[\mathbb{E}\big[(1-\rho^2)^2Z_if^{(4)}(S)+(1-\rho^2)\big((1-\rho^2)+4\rho Z_i\big)f^{(3)}(S)\nonumber\\
	&\quad+\big((\rho(\mu^2_Y+\mu_X^2)-(1+\rho^2)\mu_X\mu_Y+3\rho(1-\rho^2))+(6\rho^2-2)Z_i\big)f''(S)\nonumber\\
	&\quad+\big((2\rho\mu_X\mu_Y-\mu^2_X-\mu_Y^2+3\rho^2-1)-4\rho Z_i\big)f'(S)\,\big|\,\{Z_j\}_{j\not=i}\big]\big]\nonumber\\
 &=-\mathbb{E}\big[(1-\rho^2)^2Sf^{(4)}(S)+(1-\rho^2)\big(n(1-\rho^2)+4\rho S\big)f^{(3)}(S)\nonumber\\
	&\quad+\big(n(\rho(\mu_X^2+\mu^2_Y)-(1+\rho^2)\mu_X\mu_Y+3\rho(1-\rho^2))+(6\rho^2-2)S\big)f''(S)\nonumber\\
	&\quad+\big(n(2\rho\mu_X\mu_Y-\mu^2_X-\mu_Y^2+3\rho^2-1)-4\rho S\big)f'(S)\big],\label{nhu}
\end{align}
and substituting $f(x)=g(x/n)$ into (\ref{nhu}) yields the equation (\ref{steinc}) in the case $\sigma_X=\sigma_Y=1$. Finally, the extension to $\sigma_X,\sigma_Y>0$ follows from a simple rescaling. \qed

\vspace{2mm}

\noindent \emph{Sufficiency.} We consider the case $\sigma_X=\sigma_Y=1$; the extension to the general case follows by a straightforward rescaling. Let $W$ be a real-valued random variable such that $\mathbb{E}|W|<\infty$. On taking $f(x)=\mathrm{e}^{\mathrm{i}tx}$ in the characterising equation (\ref{steinc}) and setting $\varphi(t)=\mathbb{E}[\mathrm{e}^{\mathrm{i}tW}]$ we see that the characteristic function $\varphi(t)$ satisfies the ODE
\begin{align}
&\big(-\mathrm{i}(1-\rho^2)^2t^4-4n(1-\rho^2)\rho t^3+\mathrm{i}n^2(6\rho^2-2)t^2-4\rho tn^3-\mathrm{i}n^4\big)\varphi'(t)\nonumber\\
&\quad+\big(-\mathrm{i}n(1-\rho^2)^2t^3-n^2(\rho (\mu_X^2+\mu^2_Y)-(1+\rho^2)\mu_X\mu_Y+3\rho(1-\rho^2))t^2\nonumber\\
&\quad+\mathrm{i}n^3(2\rho \mu_X\mu_Y-\mu^2_X-\mu^2_Y+3\rho^2-1)t-n^4(\mu_X\mu_Y+\rho)\big)\varphi(t)=0. \label{ode}
\end{align}
It should be noted that  $f(x)=\mathrm{e}^{\mathrm{i}tx}$ is a complex-valued function; here we applied the characterising equation (\ref{steinc}) to the real and imaginary parts of $f$, which are real-valued functions that satisfy the conditions in the statement of the theorem.
Solving (\ref{ode}) subject to the condition $\varphi(0)=1$, we obtain
\begin{align*}
\varphi(t)=\frac{1}{([1-(1+\rho)\mathrm{i}t/n][1+(1-\rho)\mathrm{i}t/n])^{n/2}}\exp\bigg(\frac{-(\mu_X^2+\mu_Y^2-2\rho\mu_X\mu_Y)t^2/n+2\mu_X\mu_Y\mathrm{i}t}{2[1-(1+\rho)\mathrm{i}t/n][1+(1-\rho)\mathrm{i}t/n]}\bigg),
\end{align*}
which is the characteristic function of $\overline{Z}_n$ when $\sigma_X=\sigma_Y=1$ (that this formula is indeed the characteristic function of $\overline{Z}_n$ is easily seen from the formula for the characteristic function of $Z$ (see \cite[equation (10)]{craig}) and the standard formula for the characteristic function of sums of independent random variables). Thus, by the uniqueness of characteristic functions, $W=_d\overline{Z}_n$.
 \qed

\vspace{3mm}

\noindent{\emph{Proof of Theorem \ref{cor2.4}.}}  The proof of sufficiency is similar to that of Theorem \ref{cor2.2}, so is omitted. The extension from the case $n=1$ to $n\geq1$ is also similar to the argument used in the proof of Theorem \ref{cor2.2}, so is likewise omitted. We will therefore just prove the necessity part of the Stein characterisation in the case $\sigma_X=\sigma_Y=1$, with the extension to $\sigma_X,\sigma_Y>0$ following from the same argument as the one we used in proving Theorem \ref{cor2.2}. 

Suppose that $\mu_X=\mu_Y=\mu$ (recall that we are setting $\sigma_X=\sigma_Y=1$). The proof of necessity of the Stein characterisation is identical to that of Theorem \ref{cor2.2} (by setting $\mu_X=\mu_Y=\mu$) up to and including equation (\ref{F}), which we repeat here in the case $\mu_X=\mu_Y=\mu$:
\begin{align}
	&\mathbb{E}[(1-\rho^2)Zf''(Z)]+\mathbb{E}[(1-\rho^2+2\rho Z)f'(Z)]+\mathbb{E}[(\rho-Z)f(Z)]\nonumber\\
	&+\mathbb{E}[\mu(X+\mu)f(Z)]+\mathbb{E}[\mu(1-\rho)(X+\mu)f'(Z)]=0\label{Ff}.
\end{align}
By conditioning on $X$ and $V$, we have 
\begin{align}
	&\mathbb{E}[\mu(X+\mu)f(Z)]\nonumber\\
 &=\mu^2\mathbb{E}[f(Z)]+\mu\mathbb{E}\big[\mathbb{E}[Xf(Z)|V]\big]\nonumber\\
	&=\mu^2\mathbb{E}[f(Z)]+2\rho\mu\mathbb{E}[(X+\mu)f'(Z)]+\sqrt{1-\rho^2}\mu\mathbb{E}\big[\mathbb{E}[f'(Z)V|X]\big]+\mu^2(1-\rho)\mathbb{E}[f'(Z)]\nonumber\\
	&=\mu^2\mathbb{E}[f(Z)]+2\rho\mu\mathbb{E}[(X+\mu)f'(Z)]+(1-\rho^2)\mu\mathbb{E}[(X+\mu)f''(Z)]+\mu^2(1-\rho)\mathbb{E}[f'(Z)],\label{Xterm}
\end{align}
where the second equality follows by using (\ref{snsc}) with $g(x)=f(z)$ and (\ref{partial3}), and the final equality follows by using (\ref{snsc}) with $g(v)=f'(z)$ and (\ref{partial2}). On substituting (\ref{Xterm}) into (\ref{Ff}) we obtain 
\begin{align}
	&\mathbb{E}[(1-\rho^2)Zf''(Z)]+\mathbb{E}[(1-\rho^2+2\rho Z)f'(Z)]+\mathbb{E}[(\rho-Z)f(Z)]+\mu^2\mathbb{E}[f(Z)]\nonumber\\
	&+\mu^2(1-\rho)\mathbb{E}[f'(Z)]+\mu(\rho+1)\mathbb{E}[(X+\mu)f'(Z)]+(1-\rho^2)\mu\mathbb{E}[(X+\mu)f''(Z)]=0.\label{F2}
\end{align}
By multiplying (\ref{Ff}) by $(1+\rho)$ with $f$ replaced by $f'$ and subtracting (\ref{F2}), we obtain 
\begin{align*}
	&\mathbb{E}[(1-\rho^2)(1+\rho)Zf^{(3)}(Z)]+\mathbb{E}[(1+\rho)(1-\rho^2+(3\rho-1)Z)f''(Z)]\\
	&+\mathbb{E}[(2\rho^2-3\rho Z+\rho-1-Z-\mu^2+\rho\mu^2)f'(Z)]+\mathbb{E}[(Z-\mu^2-\rho)f(Z)]=0,
\end{align*}
which we recognise as the characterising equation (\ref{charac4}) in the case $\sigma_X=\sigma_Y=1$ and $n=1$.
\qed

\section*{Acknowledgements}
We would like to thank the reviewer for helpful comments. RG is funded in part by EPSRC grant EP/Y008650/1. SL is supported by a University of Manchester Research Scholar Award. HS was supported by an EPSRC PhD studentship and is currently supported by EPSRC grant EP/Y008650/1.

\end{document}